\documentclass[preprint,12pt]{elsarticle}
\usepackage{subfig}
\usepackage{wrapfig}
\usepackage{amsmath, amssymb, epsfig}
\usepackage[scaled]{helvet} 
\usepackage{fourier}
\usepackage{bm}
\usepackage{color}
\usepackage{fullpage} 	

\newcommand{\hide}[1]{}

\newcommand{\x}{\mathbf{x}}

\newcommand{\z}{\mathbf{z}}

\newcommand{\zero}{\mathbf{0}}

\newcommand{\defby}{\overset{\mathrm{\scriptscriptstyle{def}}}{=}}
\newcommand{\bigO}{\mathrm{O}}
\DeclareMathOperator*{\argmin}{arg min} 
 
\newcommand{\RR}{\mathbb{R}}
\newcommand{\NN}{\mathbb{N}}

\newcommand{\e}{\ensuremath{{\rm e}}}
\DeclareMathOperator{\EE}{\mathbb{E}}

\newcommand{\pinv}[1]{ {#1}^\dagger}
\newcommand{\norm}[1]{\ensuremath{\left\|#1\right\|_2}}

\newcommand{\rank}[1]{\ensuremath{\mathrm{\textbf{{\footnotesize rank}}}\left(#1\right)}}

\newcommand{\vct}[1]{\bm{#1}}
\newcommand{\mtx}[1]{\bm{#1}}

\newcommand{\ignore}[1]{}

\newcommand{\Id}{\mathbf{I}}

\newcommand{\zeromtx}{\mathbf{0}}

\newcommand{\mat}[1]{ {\ensuremath{\mtx{#1} }}}

\def\e{{\mathbf e}}

\def\v{{\mathbf v}}

\def\matA{\mat{A}}

\def\matI{\mat{I}}

\def\matP{\mat{P}}

\def\matU{\mat{U}}
\def\matV{\mat{V}}

\def\matSig{\mat{\Sigma}}

\newcommand\remove[1]{}

\newcommand{\vecb}{{\vct{b} }}
\newcommand{\bc}{{\vecb_{\mathcal{R}(\matA)^\bot } }}
\newcommand{\br}{{\vecb_{\mathcal{R}(\matA) } }}

\def\xls{\x_{\text{\tiny LS}}}
\def\xlsn{\overline{\x}_{\text{\tiny LS}}}

\newcommand{\colspan}[1]{\mathcal{R}(#1)}

\usepackage{palatino}

\usepackage{hyperref}

\usepackage{algorithmicx}
\usepackage[ruled]{algorithm}
\usepackage{algpseudocode}

\newcommand{\floor}[1]{\lfloor #1 \rfloor}
\newtheorem{definition}{Definition}
\newtheorem{theorem}{Theorem}
\newtheorem{proposition}[theorem]{Proposition}
\newtheorem{lemma}[theorem]{Lemma}
\newtheorem{corollary}[theorem]{Corollary}

\newtheorem{remark}{Remark}

\begin{document}
\begin{frontmatter}
\title{Randomized Block Kaczmarz Method with Projection for Solving Least Squares}

\author[deanna]{Deanna Needell\corref{cor1}}
\ead{dneedell@cmc.edu}
\author[ran]{Ran Zhao}
\ead{ran.zhao@cgu.edu}
\author[tasos]{Anastasios Zouzias}
\ead{azo@zurich.ibm.com}
\cortext[cor1]{Corresponding Author}

\address[deanna]{Dept. of Mathematical Sciences, Claremont McKenna College, Claremont, CA 91711}
\address[ran]{Dept. of Mathematics, Claremont Graduate Univ., Claremont, CA 91711}
\address[tasos]{IBM Research Lab, Zurich}

\date{\today}


\begin{abstract}
	The Kaczmarz method is an iterative method for solving overcomplete linear systems of equations $\mtx{A}\vct{x}=\vct{b}$.  The randomized version of the Kaczmarz method put forth by Strohmer and Vershynin iteratively projects onto a randomly chosen solution space given by a single row of the matrix $\mtx{A}$ and converges linearly\footnote{Mathematicians often refer to this type of convergence as exponential.} in expectation to the solution of a \textit{consistent} system.  In this paper we analyze two block versions of the method each with a randomized projection, designed to converge in expectation to the least squares solution, often faster than the standard variants.  Our approach utilizes both a row and column-paving of the matrix $\mtx{A}$ to guarantee linear convergence when the matrix has consistent row norms (called nearly standardized), and a single column-paving when the row norms are unbounded.  The proposed methods are an extension of the block Kaczmarz method analyzed by Needell and Tropp and the Randomized Extended Kaczmarz method of Zouzias and Freris.  The contribution is thus two-fold; unlike the standard Kaczmarz method, our results demonstrate convergence to the least squares solution of inconsistent systems (both methods in the nearly standardized case and the second method in other cases).  By using appropriate blocks of the matrix this convergence can be significantly accelerated, as is demonstrated by numerical experiments. 
\end{abstract}

\end{frontmatter}



%
%
%

\section{Introduction}

The Kaczmarz method~\cite{K37:Angena} is a popular iterative solver of overdetermined systems of linear equations $\mtx{A}\vct{x}=\vct{b}$.  Because of its simplicity and performance, the method and its derivatives are used in a range of applications from 
image reconstruction to digital signal processing~\cite{CFMSS92:New-Variants,Nat01:Mathematics-Computerized,SS87:Applications}.  The method performs a series of orthogonal projections and iteratively converges to the solution of the system of equations.  It is therefore computationally feasible even for very large and overdetermined systems.  

Given a vector $\vct{b}$ and an $n \times d$ full rank (real or complex) matrix $\mtx{A}$ with rows $\vct{a}_1, \vct{a}_2, \ldots \vct{a}_n$, the algorithm begins with an initial estimate $\vct{x}_0$ and cyclically projects the estimation onto each of the solution spaces.  This process can be described as follows:

$$
\vct{x}_j = \vct{x}_{j-1} + \frac{\vct{b}[i] - \langle \vct{a}_i, \vct{x}_{j-1}\rangle}{\|\vct{a}_i\|_2^2}\vct{a}_i,
$$
where $\vct{b}[i]$ denotes the $i$th coordinate of $\vct{b}$, $\vct{x}_j$ is the estimation in the $j$th iteration, $\|\cdot\|_2$ denotes the usual $\ell_2$ vector norm, $\langle \cdot, \cdot\rangle$ the standard $\ell_2$ inner product, and $i = (j\mod n) + 1$ cycles through the rows of $\mtx{A}$. 

Since the method cycles through the rows of $\mtx{A}$, the performance of the algorithm may depend heavily on the ordering of these rows.  A poor ordering may lead to very slow convergence.  To overcome this obstacle, one can select the row $\vct{a}_i$ at random to improve the convergence rate~\cite{HM93:Algebraic-Reconstruction,Nat01:Mathematics-Computerized}.  Strohmer and Vershynin proposed and analyzed a method which selects a given row with probability proportional to its $\ell_2$ norm~\cite{SV09:Randomized-Kaczmarz,SV06:Arandom}.  They show that with this selection strategy, the randomized Kaczmarz method has an expected linear convergence rate to the unique solution\footnote{Note that if the full-rank assumption is removed, the method converges to the solution set at the same linear rate with $\|\mtx{A}^{-1}\|$ replaced by the reciprocal of the smallest non-zero singular value.  See e.g. (3) of~\cite{liu2014asynchronous}.} $\vct{x}_\star$:

\begin{equation}\label{eqn:RVrate}
\mathbb{E}\|\vct{x}_j - \vct{x}_\star\|_2^2 \leq \left(1 - \frac{1}{R}\right)^j\|\vct{x}_0 - \vct{x}_\star\|_2^2,
\end{equation}
where $R$ is the scaled condition number, $R = \|\mtx{A}^{-1}\|^2\|\mtx{A}\|_F^2$, $\|\cdot\|_F$ denotes the Frobenius norm, and $\|\mtx{A}^{-1}\| \defby \inf\{M : M\|\mtx{A}\vct{x}\|_2 \geq \|\vct{x}\|_2 \text{ for all }x\}$ is well-defined since $\mtx{A}$ has full column rank.  This convergence rate~\eqref{eqn:RVrate} is essentially independent of the number of rows of $\mtx{A}$ and shows that for well-conditioned matrices, the randomized Kaczmarz method converges to the solution in just $\bigO(d)$ iterations~\cite{SV09:Randomized-Kaczmarz}.  This yields an overall runtime of $\bigO(d^2)$ which is much superior to others such as $\bigO(nd^2)$ for Gaussian elimination.  There are also cases where randomized Kaczmarz even outperforms the conjugate gradient method, see the discussion in ~\cite{SV09:Randomized-Kaczmarz} for details.

When the system is perturbed by noise or no longer consistent, $\mtx{A}\vct{x}_\star + \vct{e} = \vct{b}$, the randomized Kaczmarz method still provides expected linear convergence down to an error threshold~\cite{Nee10:Randomized-Kaczmarz}, 
\begin{equation}\label{RKnoise}
\mathbb{E}\|\vct{x}_j - \vct{x}_\star\|_2 \leq \left(1 - \frac{1}{R}\right)^{j/2}\|\vct{x}_0 - \vct{x}_\star\|_2 + \sqrt{R}\cdot\max_i\frac{|e[i]|}{\|\vct{a}_i\|_2},
\end{equation}
where $|e[i]|$ denotes the $i$th entry of $\vct{e}$. 
This result is sharp, and shows that the randomized Kaczmarz method converges with a radius proportional to the magnitude of the largest entry of the noise in the system.  Since the iterates of the Kaczmarz method always lie in a single solution space, the method clearly will not converge to the least squares solution of an inconsistent system.

\subsection{Randomized Extended Kaczmarz}
The bound~\eqref{RKnoise} demonstrates that the randomized Kaczmarz method performs well when the noise in inconsistent systems is small.  Zouzias and Freris introduced a variant of the method which utilizes a random projection to iteratively reduce the norm of the error~\cite{ZF12:REK}.  They show that the estimate of this \textit{Randomized Extended Kaczmarz} (REK) method converges linearly in expectation to the least squares solution of the system, breaking the radius barrier of the standard method.  The algorithm maintains not only an estimate $\vct{x}_j$ to the solution but also an approximation $\vct{z}_j$ to the projection of $\vct{b}$ onto the range of $\mtx{A}$:

\begin{equation}\label{eq:rek}
\vct{x}_j = \vct{x}_{j-1} + \frac{\vct{b}[i] - \vct{z}_{j-1}[i] - \langle \vct{a}_i, \vct{x}_{j-1}\rangle}{\|\vct{a}_i\|_2^2}\vct{a}_i, \quad 
\vct{z}_j = \vct{z}_{j-1} - \frac{\langle \vct{a}^{(k)}, \vct{z}_{j-1}\rangle}{\|\vct{a}^{(k)}\|_2^2}\vct{a}^{(k)},
\end{equation}
where in iteration $j$, $\vct{a}_i$ and $\vct{a}^{(k)}$ is the row and column of $\mtx{A}$, respectively, each chosen randomly with probability proportional to their Euclidean norms.  In this setting, we no longer require that the matrix $\mtx{A}$ be full rank, and ask for the least squares solution,
$$
\vct{x}_{LS} \defby \argmin_{\vct{x}} \|\vct{b} - \mtx{A}\vct{x}\|_2 = \mtx{A}^{\dagger}\vct{b},
$$
where $\mtx{A}^{\dagger}$ denotes the Moore-Penrose pseudoinverse of $\mtx{A}$.  Zouzias and Freris showed that the REK method converges linearly in expectation to the least squares solution~\cite{ZF12:REK},
\begin{equation}\label{eqn:REKrate}
\mathbb{E}\|\vct{x}_j - \vct{x}_{LS}\|_2^2 \leq \left(1 - \frac{1}{K^2(\mtx{A})}\right)^{j/2}\left(\|\vct{x}_{LS}\|_2^2 + \frac{2\|\vct{b}\|_2^2}{\sigma_{\min}^2(\mtx{A})}\right),
\end{equation}
where $\sigma_{\min}(\mtx{A})$ is the smallest non-zero singular value of $\mtx{A}$ and $K(\mtx{A}) = \frac{\|\mtx{A}\|_F}{\sigma_{\min}(\mtx{A})}$ denotes its scaled condition number.

\subsection{The block Kaczmarz method}
Recently, Needell and Tropp analyzed a block version of the simple randomized Kaczmarz method~\cite{needell2013paved}.  
Like the traditional method, this version iteratively projects the current estimation onto the solution spaces.  However, rather than using the solution space of a \textit{single} equation, the block method projects onto the solution space of many equations simultaneously by selecting a block of rows rather than a single row.  For a subset $\tau \subset \{1, 2, \ldots, n\}$, denote by $\mtx{A}_\tau$ the submatrix of $\mtx{A}$ whose rows are indexed by $\tau$.  We again begin with an arbitrary guess $\vct{x}_0$ for the solution of the system.  Then for each iteration $j\geq 1$, select a block $\tau = \tau_j$ of rows.  To obtain the next iterate, we project the current estimation onto the solution space of the equations listed in $\tau$~\cite{needell2013paved}:

\begin{equation}\label{eq:iterates}
\vct{x}_j = \vct{x}_{j-1} + (\mtx{A_\tau})^\dagger(\vct{b}_\tau - \mtx{A}_\tau\vct{x}_{j-1}).
\end{equation}

Here, the conditioning of the blocks $\mtx{A_\tau}$ plays a crucial row in the behavior of the method.  Indeed, if each block is well-conditioned, its pseudoinverse can be applied efficiently using an iterative method such as CGLS~\cite{Bjo96:Numerical-Methods}.  To guarantee such properties, Needell and Tropp utilize a \textit{paving} of the matrix $\mtx{A}$.

\begin{definition}[Row Paving]\label{def:paving}
A $(p, \alpha, \beta)$ \emph{row paving} of a $n\times d$ matrix $\mtx{A}$ is a partition $\mathcal{T} = \{ \tau_1, \dots, \tau_p \}$ of the rows such that
$$
\alpha \leq \lambda_{\min}( \matA_\tau \matA_{\tau}^* )
\quad\text{and}\quad
\lambda_{\max}( \matA_{\tau} \matA_{\tau}^* ) \leq \beta
\quad\text{for each $\tau \in T$,}
$$
where again we denote by $\matA_{\tau}$ the $ |\tau|\times d$ submatrix of $\matA$.  We refer to the number $p$ as the \emph{size} of the paving, and the numbers $\alpha$ and $\beta$ are called the \emph{lower} and \emph{upper paving bounds}.  
\end{definition} 

We refer to a row paving of $\matA^*$ as a \textit{column paving} of $\matA$.  We thus seek pavings of $\matA$ with small number of blocks $p$ and upper paving constant $\beta$.  In the following, we will assume one has access to such a paving, and discuss in Section~\ref{sec:pavings} how to construct such pavings for various types of matrices.  When the matrix has unit-norm rows, equipped with such a row paving of $\mtx{A}$, the main result of~\cite{needell2013paved}
shows that the randomized block Kaczmarz algorithm~\eqref{eq:iterates} exhibits linear convergence in expectation:  

\begin{equation}\label{blockRate}
\mathbb{E}\|\vct{x}_j - \vct{x}_\star\|_2^2 \leq \left(1 - \frac{1}{C'\kappa^2(\mtx{A})\log (1+n)}\right)^j \|\vct{x}_0 - \vct{x}_\star\|_2^2 + \frac{3\|\vct{e}\|_2^2}{\sigma^2_{\min}(\mtx{A})},
\end{equation}
where $C'$ is an absolute constant and $\kappa(\mtx{A}) = \frac{\sigma_{\max}(\mtx{A})}{\sigma_{\min}(\mtx{A})}$ is the condition number of $\mtx{A}$.

Since each iteration of the block method utilizes multiple rows, one can compare the rate of~\eqref{blockRate} and~\eqref{eqn:RVrate} by considering convergence per epoch (one cycle through the rows of $\mtx{A}$).  From this analysis, one finds the bounds to be comparable.  However, the block method can utilize fast matrix multiplies and efficient implementation, yielding dramatic improvements in computational time.  See~\cite{needell2013paved} for details and empirical results.

\subsection{Contribution}
The REK method breaks the so-called \textit{convergence horizon} of standard Kaczmarz method, allowing convergence to the least squares solution of inconsistent systems.  The block Kaczmarz method on the other hand, allows for significant computational speedup and accelerated convergence to within a fixed radius of the least squares solution.  
The main contribution of this paper analyzes a randomized block Kaczmarz method which also incorporates a blocked projection step, which provides accelerated convergence to the least squares solution.  In this case we need a column partition for the projection step and a row partition for the Kaczmarz step.  
Our results show that this method offers both linear convergence to the least squares solution $\xls$ when a row paving can be obtained, and improved convergence speed due to the blocking of both the rows and columns.  In addition, we present a block coordinate descent variant which utilizes only a column paving and also yields linear convergence.  We will see that the desired column paving may be obtained easily for arbitrary matrices, so this variant is especially useful when a good row paving cannot be obtained.

\subsection{Organization}
In Section~\ref{sec:doublerbk} we begin by presenting the double block randomized Kaczmarz method which utilizes both a row and column paving, and provides the strongest theoretical results overall.  We discuss methods for obtaining the desired matrix pavings in Section~\ref{sec:pavings}. Section~\ref{sec:coord} introduces the variant of the block and extended Kaczmarz methods that requires only a column paving of the matrix, readily accessible for arbitrary matrices.  In Section~\ref{sec:exper} we present some experimental results for the various algorithms.  We conclude with a discussion of related work and open directions in Section~\ref{sec:discuss}.  The appendix includes proofs of intermediate results used along the way.

\section{The Randomized Double Block Kaczmarz Method}\label{sec:doublerbk}

It is natural to ask whether one can consider both a row partition and a column partition in the Kaczmarz method, blocking both in the Kaczmarz update step and the projection step.  Indeed, utilizing blocking in both steps yields Algorithm~\ref{alg:LSrek} below. 
We thus propose the following randomized block extended Kaczmarz method using double partitioning.  We will see later that since this method requires a row paving, it works very well when each of the row norms of the matrix are relatively similar (see Sections \ref{sec:pavings} and \ref{sec:sum}).

\begin{algorithm}{}
	\caption{Randomized Double Block Kaczmarz Least Squares Solver }\label{alg:LSrek}
\begin{algorithmic}[1]
\Procedure{}{$\matA$, $\vecb$, $T$, $\mathcal{T}$, $\mathcal{S}$}\Comment{$\matA\in\RR^{n\times d}$, $\vecb\in\RR^n$, $T\in \NN$, column partition $\mathcal{T}$ of $[d]$, row partition $\mathcal{S}$ of $[n]$}
\State Initialize $\x_0 = \zeromtx$ and $\z_0 =\vecb$
\For {$k=1,2,\ldots, T  $ }
	\State Pick $\tau_k\in\mathcal{T}$ and $\upsilon_k\in\mathcal{S}$ uniformly at random 
	\State Set $ \z_{k} = \z_{k-1} - {\matA}_{\tau_k}\pinv{({\matA}_{\tau_k})} \z_{k-1}$  \Comment{${\matA}_{\tau_k}$: $n\times |\tau_k|$ submatrix of ${\matA}$}
	\State Update $\vct{x}_k = \vct{x}_{k-1} + (\mtx{A}_{\upsilon_k})^\dagger(\vct{b}_{\upsilon_k} - (\vct{z}_{k})_{\upsilon_k} - \mtx{A}_{\upsilon_k}\vct{x}_{k-1})$  \Comment{$\matA_{\upsilon_k}$: $|\upsilon_k| \times d$ submatrix of $\matA$}
\EndFor
\State Output $\x_T$
\EndProcedure
\end{algorithmic}
\end{algorithm}

Combining the theoretical approaches in~\cite{needell2013paved,ZF12:REK} we will prove the following result about the convergence of Algorithm~\ref{alg:LSrek}.  This result utilizes both a column paving and row paving.

\begin{theorem}\label{thm:LSrek}
Algorithm~\ref{alg:LSrek} with input $\matA$, $\vecb$, $T\in\NN$, $(\overline{p}, \overline{\alpha}, \overline{\beta})$ column paving $\mathcal{T}$ of ${\matA}$, and $(p, \alpha, \beta)$ row paving $\mathcal{S}$ of $\matA$, outputs an estimate vector $\x_{T}$ that satisfies
\[ 
\EE \|\vct{x}_T - \xls\|_2^2 \leq \gamma^T\|\vct{x}_0 - \xls\|_2^2 + \left({\gamma}^{\floor{T/2}} + \overline{\gamma}^{\floor{T/2}}\right)\frac{\norm{\br}^2}{\alpha(1-\gamma)},
\]
where $\gamma = 1 -\frac{\sigma_{\min}^2({\matA})}{ {p} {\beta} }$ and $\overline{\gamma} = 1 -\frac{\sigma_{\min}^2({\matA})}{ \overline{p} \overline{\beta} }$.
\end{theorem}

\textit{Proof. }
We begin with the following lemma which is motivated by Lemma~$2.2$ of~\cite{needell2013paved}, and shows that the iterates $\z_k$ converge linearly to the projection of $\vecb$ onto the kernel of $\matA^*$.

\begin{lemma}\label{thm:randOP:block}
Let $\vecb$ be a fixed vector and $\mathcal{T}$ be a $(\overline{p}, \overline{\alpha}, \overline{\beta})$ column paving of $\matA$.  Assuming the notation of Theorem~\ref{thm:LSrek}, for every $k>0$ it holds that
\begin{equation}\label{eq:main}
\EE \norm{\z_k - \bc }^2 \leq \left(1 -\frac{\sigma_{\min}^2(\matA)}{ \overline{p}\overline{\beta} }\right)^k  \norm{\br}^2.
\end{equation}
where $\bc := (\Id - \matA\pinv{\matA}) \vecb$.
\end{lemma}
\textit{Proof. }
Let $\matP_{\tau_k} = \matA_{\tau_k}\pinv{(\matA_{\tau_k})}$ and notice $\z_k = (\matI-\matP_{\tau_k}) \z_{k-1}$. Define $\e_k = \z_k - \bc$ for $k\geq 0$.  Then,
\begin{align*}
	\e_k & = (\matI - \matP_{\tau_k}) \z_{k-1} - \bc \\
	& = (\matI - \matP_{\tau_k}) \z_{k-1} - (\matI - \matP_{\tau_k}) \bc \\
	& = (\matI - \matP_{\tau_k}) \e_{k-1},
\end{align*}
where the first equality follows by the definition of $\z_k$, the second by orthogonality between the range of $\matP_{\tau_k}$ and $\bc$, and the final equality by definition of $\e_{k-1}$. Next, we prove that 
\[\EE_{k-1} \norm{\e_k}^2 \leq \left(1 -\frac{\sigma_{\min}^2(\matA)}{ \overline{p}\overline{\beta} }\right) \norm{\e_{k-1}}^2\]
where $\EE_{k-1}$ is the expectation conditioned over the first $(k-1)$ iterations of the algorithm. By orthogonality of the projector $\matP_{\tau_k}$ and the Pythagorean theorem, one has $\norm{(\matI - \matP_{\tau_k}) \e_{k-1}}^2 = \norm{\e_{k-1}}^2 - \norm{\matP_{\tau_k} \e_{k-1}}^2$, hence it suffices to lower bound $\EE_{k-1} \norm{\matP_{\tau_k} \e_{k-1}}^2$. Let $\matA_{\tau_k} := \matU_{\tau_k} \matSig_{\tau_k} \matV_{\tau_k}^* $ be the truncated SVD decomposition of $\matA_{\tau_k}$ where $\matSig_{\tau_k}$ is a $\rank{\matA_{\tau_k}} \times \rank{\matA_{\tau_k}}$ diagonal matrix containing the non-zero singular values of $\matA_{\tau_k}$. Then,
\begin{align*}
	\EE_{k-1} \norm{\matP_{\tau_k} \e_{k-1}}^2 & = \EE_{k-1} \norm{ \matU_{\tau_k}\matU_{\tau_k}^* \e_{k-1}}^2 \\
	& = \EE_{k-1} \norm{ \matU_{\tau_k}^* \e_{k-1}}^2 \\
											   & = \EE_{k-1} \norm{\matSig_{\tau_k}^{-1} \matV_{\tau_k}^* \matA_{\tau_k}^* \e_{k-1}}^2 \\
											   & \geq \EE_{k-1} \sigma^2_{\min}(\matSig_{\tau_k}^{-1}\matV_{\tau_k}^*)\norm{ \matA_{\tau_k}^* \e_{k-1}}^2 \\
										       & = \EE_{k-1} \frac{\norm{ \matA_{\tau_k}^* \e_{k-1}}^2}{\sigma^2_{\max}( \matSig_{\tau_k})} \\
										       & \geq \frac1{\overline{p}\overline{\beta}}\sum_{\tau_k\in\mathcal{T}}\norm{ \matA_{\tau_k}^* \e_{k-1}}^2 \\
										       & = \frac1{\overline{p}\overline{\beta}}\norm{ \matA^* \e_{k-1} }^2 \\
										       & \geq \frac{\sigma^2_{\min}(\matA)}{\overline{p}\overline{\beta}}\norm{\e_{k-1}}^2,
\end{align*}
where the first equality follows since $\matP_{\tau_k} = \matU_{\tau_k} \matU_{\tau_k}^*$, the second equality from the unitary invariance property of the 
Euclidean norm, the third equality by replacing $\matU_{\tau_k}^*$ with $\matSig_{\tau_k}^{-1}\matV_{\tau_k}^* \matA_{\tau_k}^* $, the next three lines follow 
since $\sigma_{\min}^2(\matSig_{\tau_k}^{-1} \matV_{\tau_k}^*) = 1/\sigma_{\max}^2(\matSig_{\tau_k}) = 1/ \sigma^2_{\max}(\matA_{\tau_k})$ and by the paving assumption, and the final inequality follows since $\e_{k}\in \colspan{\matA}$ for all $k\geq 0 $ (indeed, $\e_0 = \br\in\colspan{\matA}$ and it follows that $\e_k \in \colspan{\matA}$ for every $k\geq 0$ by the recursive definition of $\e_k$). It follows that 
\begin{equation}
	\EE_{k-1} \norm{ \e_{k}}^2  = \norm{\e_{k-1} }^2 - \EE_{k-1} \norm{\matP_{\tau_k} \e_{k-1}}^2 \leq \left(1 -\frac{\sigma_{\min}^2(\matA)}{ \overline{p}\overline{\beta} }\right) \norm{\e_{k-1}}^2.
\end{equation}
Repeat the above inequality $k$ times and notice that $\e_0 = \vct{b} - \bc = \br$ to conclude.
$\hfill\Box$

Next, Lemma 2.2 of~\cite{needell2013paved} shows that for any vector $\vct{u}$,
\begin{equation}\label{eq:NT}
\EE \norm{\left(\Id - (\mtx{A}_{\tau})^\dagger\mtx{A}_{\tau}\right)\vct{u}}^2 \leq \left( 1 - \frac{\sigma^2_{\min}(\matA)}{p\beta}\right)\|\vct{u}\|_2^2.
\end{equation}

Since the range of $\Id - (\mtx{A}_{\tau})^\dagger\mtx{A}_{\tau}$ and $(\mtx{A}_{\tau})^\dagger$ are orthogonal, we have
\begin{equation}\label{eq:orth}
\|\vct{x}_k - \xls\|_2^2 = \norm{\left(\Id - (\mtx{A}_{\tau})^\dagger\mtx{A}_{\tau}\right)(\vct{x}_{k-1} - \xls)}^2 + \norm{(\mtx{A}_{\tau})^\dagger((\vct{z}_{k})_{\tau} - \vct{b}^{\perp}_{\tau})}^2,
\end{equation}
where for shorthand we will write $\vct{b}^{\perp}_{\tau}$ to mean $(\bc)_{\tau}$.
Combining~\eqref{eq:NT} with $\vct{u} = \vct{x}_{k-1} - \xls$ along with~\eqref{eq:orth}, we have
\begin{align}\label{eq:bnd1}
\EE \|\vct{x}_k - \xls\|_2^2 &\leq \left( 1 - \frac{\sigma^2_{\min}(\matA)}{p\beta}\right)\EE\norm{\vct{x}_{k-1} - \xls}^2 + \EE\norm{(\mtx{A}_{\tau})^\dagger((\vct{z}_{k})_{\tau} - \vct{b}^{\perp}_{\tau})}^2\notag\\
&\leq \left( 1 - \frac{\sigma^2_{\min}(\matA)}{p\beta}\right)\EE\norm{\vct{x}_{k-1} - \xls}^2 + \frac{1}{\sigma^2_{\min}(\matA_{\tau})}\EE\norm{(\vct{z}_{k})_{\tau} - \vct{b}^{\perp}_{\tau}}^2\notag\\
&\leq \left( 1 - \frac{\sigma^2_{\min}(\matA)}{p\beta}\right)\EE\norm{\vct{x}_{k-1} - \xls}^2 + \frac{1}{\alpha}\EE\norm{\vct{z}_{k} - \bc}^2.
\end{align}

To apply this bound recursively, we will utilize an elementary lemma.  It is essentially proved in~\cite[Theorem 8]{ZF12:REK} but for completeness we recall its proof in the appendix. 

\begin{lemma}\label{lem:elem}
Suppose that for some $\gamma, \overline{\gamma} < 1$, the following bounds hold for all $k^* \geq 0$:
\begin{equation}\label{eq:bnds}
\EE \|\vct{x}_{k^*} - \xls\|_2^2 \leq \gamma\EE\|\vct{x}_{k^*-1} - \xls\|_2^2 + r_{k^*} \quad\text{and}\quad r_{k^*} \leq \overline{\gamma}^{k^*}B.
\end{equation}
Then for any $T > 0$,
$$
\EE \|\vct{x}_T - \xls\|_2^2 \leq \gamma^T\|\vct{x}_0 - \xls\|_2^2 + \left(\gamma^{\floor{T/2}} + \overline{\gamma}^{\floor{T/2}}\right)\frac{B}{1-\gamma}.
$$

\end{lemma}

Using $r_k = \frac{1}{\alpha}\EE\norm{\vct{z}_{k} - \bc}^2$, $B = \frac{\norm{\br}^2}{\alpha}$, $\gamma = 1 -\frac{\sigma_{\min}^2({\matA})}{ {p} {\beta} }$, and $\overline{\gamma} = 1 -\frac{\sigma_{\min}^2({\matA})}{ \overline{p} \overline{\beta} }$, we see by Lemma~\ref{thm:randOP:block} and~\eqref{eq:bnd1} that the bounds~\eqref{eq:bnds} hold.  Applying Lemma~\ref{lem:elem} completes the proof.

$\hfill\Box$

\subsection{Implementation}
In Section~\ref{sec:pavings} we will see that matrix row-pavings for \textit{standardized} matrices, those whose rows have unit norm, can be obtained readily.  One can thus use the column-normalized version of the matrix $\matA$ in line 5 of Algorithm~\ref{alg:LSrek} and the corresponding column paving.  Note that one need not have access to the complete column-standardized matrix, instead the columns of the submatrix in line 5 could be normalized on the fly.   See Section~\ref{sec:pavings} for more on obtaining row pavings and further details. 

\subsection{Comparison of Convergence Rates}

The bound of Theorem~\ref{thm:LSrek} improves upon that of the randomized block Kaczmarz method because it demonstrates linear convergence to the least squares solution $\xls$, whereas \eqref{blockRate} shows convergence only within a radius proportional to $\|\e\|_2^2$, which we call the \textit{convergence horizon}.  Algorithm~\ref{alg:LSrek} is able to break this barrier because it iteratively removes the component of $\vecb$ which is orthogonal to the range of $\matA$.  This of course is also true of the randomized Extended Kaczmarz method~\eqref{eq:rek} as it also breaks this horizon barrier. To compare the rate of~\eqref{eqn:REKrate} to that of Theorem~\ref{cor:LSrek}, we consider two important scenarios. 

First, consider the case when $\matA$ is nearly square, and each submatrix can be applied efficiently via a fast multiply.  In this case, each iteration of Algorithm~\ref{alg:LSrek} incurs approximately the same computational cost as an iteration of the REK method.  Thus, we may directly compare the convergence rates of Theorem~\ref{cor:LSrek} and~\eqref{eqn:REKrate} to find that Algorithm~\ref{alg:LSrek} is about $n/(p\beta)$ times faster than REK in this setting.  Thus when $n$ is much larger than $p\beta$, this can result in a significant speedup.

Alternatively, if the matrix $\matA$ does not admit a fast multiply, it is fair to only compare the convergence rate per epoch, since each iteration of Algorithm~\ref{alg:LSrek} may require more computational cost than those of REK.  Since an epoch of Algorithm~\ref{alg:LSrek} and REK consist of $p$ and $n$ iterations, respectively, we see that the rate of the former is proportional to $\sigma^2_{\min}(\matA)/\beta$ whereas that of REK is proportional to $\sigma^2_{\min}(\matA)$.  We see in this case that these bounds suggest REK exhibits faster convergence (assuming $\beta > 1$).  However, as observed in the randomized Block Kaczmarz method, the block methods still display faster convergence than their single counterparts because of implicit computational issues in the linear algebraic subroutines.  See the discussion in~\cite{needell2013paved} and the experimental results below for further details.

\section{Obtaining matrix pavings}\label{sec:pavings}

We devote this section to a brief discussion about matrix pavings and how they may be obtained to utilize the results of Theorem~\ref{thm:LSrek}.  The results discussed here on matrix pavings stem from \textit{subset selection}, the problem of selecting a large submatrix with desired geoemtric properties, whose origins come from the well-known Restricted Invertibility Principle of Bourgain Tzafriri~\cite{BT87:Invertibility-Large}.  The literature now contains several results on subset selection, see e.g.~\cite{BT91:Problem-Kadison-Singer,KT94:Some-Remarks,Ver01:Johns-Decompositions,Ver06:Random-Sets,tropp2009column} and \cite{Sri10:Spectral-Sparsification,Nao11:Sparse-Quadratic,SS12:Elementary-Proof,You12:Restricted-Invertibility} for recent advancements.  We summarize here a few results useful for our purposes. 

The simplest type of matrix to first consider is one which has unit-norm rows, which we call $\textit{row-standardized}$ (and a matrix with unit-norm columns is \textit{column-standardized}).  A surprising result shows that every row-standardized matrix admits a row paving with well-controlled paving parameters.  Tropp proves the following result in~\cite[Thm.~1.2]{Tro09:Column-Subset}, whose origins are due to Bourgain and Tzafriri~\cite{BT87:Invertibility-Large,BT91:Problem-Kadison-Singer} and Vershynin~\cite{Ver06:Random-Sets}. 

\begin{proposition}[Existence of Good Pavings] \label{prop:intro-paving}
Fix a number $\delta \in (0, 1)$ and row-standardized matrix $\matA$ with $n$ rows.
Then $\matA$ admits a $(p, \alpha, \beta)$ row paving with
\begin{equation}\label{eq:paving}
p \leq C \cdot \delta^{-2} \|\matA\|^2 \log(1+n)
\quad\text{and}\quad
1 - \delta \leq \alpha \leq \beta \leq 1 + \delta,
\end{equation}
where $C$ denotes an absolute constant.
\end{proposition}

Proposition~\ref{prop:intro-paving} shows the \textit{existence} of such a paving, but the literature provides various efficient mechanisms for the construction of good pavings as well.  In many cases one constructs such a paving simply by choosing a partition of an appropriate size \textit{at random}, see~\cite{Tro09:Column-Subset} for an efficient method to compute a paving satisfying~\eqref{eq:paving}.  See also~\cite{needell2013paved} and the references therein for a thorough discussion of these types of results.

If the matrix $\matA$ is row-standardized, one can thus construct a row paving satisfying~\eqref{eq:paving}.  If the matrix also naturally admits a column paving (for example if it is symmetric or positive semi-definite), then both pavings will have such bounded parameters.  If the latter does not hold, one can instead utilize the column-standardized version of $\matA$, which we denote $\overline{\matA}$, in line 5 of Algorithm~\ref{alg:LSrek}.  The standardization can either be done on the fly during the algorithm, or ahead of time.  Either way, Proposition~\ref{prop:intro-paving} can be combined with Theorem~\ref{thm:LSrek} to yeild the following corollary.

\begin{corollary}\label{cor:LSrek}
Suppose Algorithm~\ref{alg:LSrek} is run on a row-standardized matrix $\matA$, $\vecb$, $T\in\NN$, $(\overline{p}, \overline{\alpha}, \overline{\beta})$ column paving $\mathcal{T}$ of the column-standardized version $\overline{\matA}$, and $(p, \alpha, \beta)$ row paving $\mathcal{S}$ of $\matA$, both guaranteed by Proposition~\ref{prop:intro-paving}.  Then the estimate vector $\x_{T}$ satisfies
\[ 
\EE \|\vct{x}_T - \xls\|_2^2 \leq \gamma^T\|\vct{x}_0 - \xls\|_2^2 + \left({\gamma}^{\floor{T/2}} + \overline{\gamma}^{\floor{T/2}}\right)\frac{C\norm{\br}^2}{(1-\gamma)},
\]
where $\gamma = 1 -\frac{C}{ \kappa^2(\matA)\log(1+n) }$, $\overline{\gamma} = 1 -\frac{C}{ \kappa^2(\overline{\matA})\log(1+d) }$and $\kappa(\mtx{A}) = \frac{\sigma_{\max}(\mtx{A})}{\sigma_{\min}(\mtx{A})}$ and $\kappa(\overline{\matA}) = \frac{\sigma_{\max}(\overline{\matA})}{\sigma_{\min}(\overline{\matA})}$ denote the condition numbers of $\mtx{A}$ and $\overline{\matA}$, respectively.
\end{corollary}

If the matrix $\matA$ is not row-standardized, one can run Algorithm~\ref{alg:LSrek} on the standardized system, along with a paving guaranteed by Proposition~\ref{prop:intro-paving}.  Clearly, the method will converge to the new least squares solution, which does not necessarily coincide with the original least squares solution in the inconsistent case.  

Alternatively, if the matrix $\matA$ is not row-standardized and one wishes to ensure convergence to the true least-squares solution $\xls$, Proposition~\ref{prop:intro-paving} can still be used with sub-optimal paving bounds.  Indeed, write $\tilde{\mtx{D}}$ as the diagonal matrix whose entries correspond to the reciprocals of the row-norms $\|\vct{a}_i\|_2$ of $\matA$, so that $\tilde{\matA} = \tilde{\mtx{D}}\matA$ has unit-norm rows.  Set
$$
a_{\text{max}} = \max_i \|\vct{a}_i\|_2^2 \quad\text{and}\quad a_{\text{min}} = \min_i \|\vct{a}_i\|_2^2
$$
  Utilize Proposition~\ref{prop:intro-paving} on $\tilde{\matA}$ to obtain a row-paving with parameters $\tilde{p}$, $\tilde{\alpha}$ and $\tilde{\beta}$ so that 
$$
\tilde{p} \leq C \cdot \delta^{-2} \|\tilde{\matA}\|^2 \log(1+n)
\quad\text{and}\quad
1 - \delta \leq \tilde{\alpha} \leq \tilde{\beta} \leq 1 + \delta.
$$
If one uses this same paving for $\matA$, one has (quite pessimistically) that the corresponding paving parameters $p$, $\alpha$, and $\beta$ for $\matA$ satisfy
$$
p \leq C \cdot \delta^{-2}a_{\text{min}}^{-1} \|{\matA}\|^2 \log(1+n)
\quad\text{and}\quad
a_{\text{min}}(1 - \delta) \leq {\alpha} \leq {\beta} \leq a_{\text{max}}(1 + \delta).
$$
One can then directly apply Theorem~\ref{thm:LSrek} with this paving to obtain an analogous version of Corollary~\ref{cor:LSrek} for arbitrary matrices, which shows that 
\begin{equation}\label{eq:dr}
\EE \|\vct{x}_T - \xls\|_2^2 \leq \gamma^T\|\vct{x}_0 - \xls\|_2^2 + \left({\gamma}^{\floor{T/2}} + \overline{\gamma}^{\floor{T/2}}\right)\frac{C\norm{\br}^2}{a_{\text{min}}(1-\gamma)},
\end{equation}
where now $\gamma = 1 -\frac{Ca_{\text{min}}}{ a_{\text{max}}\kappa^2(\matA)\log(1+n) } = 1 -\frac{C}{ a_r\kappa^2(\matA)\log(1+n) }$ where we have written $a_r = a_{\text{max}}/a_{\text{min}}$ to denote the dynamic range of the row norms of $\matA$ (and $\overline{\gamma}$ remains the same as in the corollary).  This demonstrates that if the dynamic range is bounded, the convergence rate is the same as in the standardized case, up to constants.  In either case, these results demonstrate linear convergence to the true least squares solution. 

An alternative to these types of bounds can be obtained by simply only using a column paving for the matrix.  As we elaborate below, the least squares solution can still be obtained from the column-normalized system.  Since column-pavings can be easily attained in this case via Proposition~\ref{prop:intro-paving}, such a method offers a nice alternative in situations where the dynamic range of the row norms is unknown or unbounded.  We propose such a method in the next section.

\section{A Randomized Block Coordinate Descent Method}\label{sec:coord}
We next present a simple variant of the extended Kaczmarz method which utilizes only a column paving of the matrix.  For that reason, one need not worry about whether the matrix $\matA$ is row-standardized.  Moreover, the column-standardized version can be used within the algorithm which guarantees bounded paving parameters, while still finding the true least squares solution.

 Utilizing the benefits of both the block variant and the randomized extension, we propose the following randomized block coordinate descent method for the inconsistent case. 

\begin{algorithm}{}
	\caption{Randomized Block Least Squares Solver }\label{alg:LSsimple}
\begin{algorithmic}[1]
\Procedure{}{$\matA$, $\vecb$, $T$, $\mathcal{T}$}\Comment{$\matA\in\RR^{n\times d}, \vecb\in\RR^n$, $T\in \NN$, column partition $\mathcal{T}$ of $[d]$}
\State Initialize $\x_0 = \zeromtx$ and $\z_0 =\vecb$
\For {$k=1,2,\ldots, T $ }
	\State Pick $\tau_k\in\mathcal{T}$ uniformly at random 
    \State Compute $\vct{w}_k = \pinv{(\matA_{\tau_k})} \z_{k-1}$  \Comment{$\matA_{\tau_k}$: $n\times |\tau_k|$ submatrix of $\matA$}
    \State Update $ (\x_{k})_{\tau_k} = (\x_{k-1})_{\tau_k} + \vct{w}_{k}$
	\State Set $ \z_{k} = \z_{k-1} - \matA_{\tau_k} \vct{w}_k$ 
\EndFor
\State Output $\x_T$
\EndProcedure
\end{algorithmic}
\end{algorithm}

\subsection{Analysis of Randomized Block Least Squares Solver}
We may utilize some of the previous analysis to prove convergence of Algorithm~\ref{alg:LSsimple}.  
Observe that Step 6 of Algorithm~\ref{alg:LSrek} is identical to Steps 5 and 7 of Algorithm~\ref{alg:LSsimple}, therefore Lemma~\ref{thm:randOP:block} implies that 
\begin{equation}\label{eq:zous2}
	\EE \norm{\z_k - \bc }^2 \leq \left(1 -\frac{\sigma_{\min}^2({\matA})}{ \overline{p} {\overline{\beta}} }\right)^k  \norm{\br}^2.
\end{equation}

To utilize this result, we aim to relate the iterates $\z_k$ to the estimation $\x_k$.  The following claim quantifies precisely this relation.
\begin{lemma}\label{claim:1}
For every $k\geq 0$, at the end of the $k$-th iteration, it holds that $\z_{k+1} = \vecb - \matA \x_{k+1}$.
\end{lemma}
\textit{Proof. }
We prove by induction on $k$. Instate the notation of Algorithm~\ref{alg:LSsimple}, and for two sets $S_1$ and $S_2$ write $S_1\setminus S_2 = S_1\cap S_2^c$ to denote set subtraction.  For the base case of $k=0$, we have $(\x_1)_{\tau_0}=\vct{w}_0$ and $(\x_1)_{[n]\setminus \tau_0} = \zero$ and moreover, $\z_1 = \z_0 - \matA_{\tau_0} \vct{w}_0= \vecb - \matA \x_1$. Assume that $\z_\ell = \vecb - \matA \x_\ell$ is true for some $\ell>0$, we will show that it holds for $\ell+1$. For the sake of notation, denote $\matP_{\ell} = \matA_{\tau_\ell}\pinv{\matA_{\tau_\ell}}$. Then
\begin{equation}\label{eq:1}
\z_{\ell+1}   =  \z_{\ell} -  \matP_{\ell}\z_{\ell}  =  \vecb - \matA \x_{\ell} -  \matP_{\ell} \z_{\ell}
\end{equation}
the first equality follows by the definition of $\z_{\ell+1}$, the second equality follows by induction hypothesis.
Now, it follows that 
\begin{eqnarray*}
	\matA \x_{\ell+1} & = & \matA_{\tau_\ell}(\x_{\ell+1})_{\tau_\ell}  + \matA_{[n]\setminus \tau_\ell} (\x_{\ell+1})_{[n]\setminus \tau_\ell} \\
					 & = & \matA_{\tau_\ell}(\x_{\ell})_{\tau_\ell}  + \matA_{\tau_\ell} \vct{w}_{\ell} + \matA_{[n]\setminus \tau_\ell} (\x_{\ell+1})_{[n]\setminus \tau_\ell} \\
					 & = & \matA_{\tau_\ell}(\x_{\ell})_{\tau_\ell}  + \matA_{\tau_\ell} \vct{w}_{\ell} + \matA_{[n]\setminus \tau_\ell} (\x_{\ell})_{[n]\setminus \tau_\ell} \\
					 & = & \matA \x_{\ell}  + \matA_{\tau_\ell} \vct{w}_{\ell}.
\end{eqnarray*}
the first equality follows by Step $6$ of the algorithm (update on $\x$), the second equality because $(\x_{\ell+1})_{[n]\setminus \tau_\ell} = (\x_{\ell})_{[n]\setminus \tau_\ell}$. Hence, $\matA \x_{\ell} = \matA \x_{\ell+1} - \matA_{\tau_\ell} \vct{w}_{\ell}$.
Now, the right hand side of~\eqref{eq:1} can be rewritten as 
\begin{eqnarray*}
	\vecb - \matA \x_{\ell} -  \matP_{\ell} \z_{\ell} & = & \vecb - \matA \x_{\ell+1} + \matA_{\tau_\ell}\vct{w}_{\ell}  -  \matP_{\ell} \z_{\ell} \\
												 & = &  \vecb - \matA \x_{\ell+1}.
\end{eqnarray*}
the last equality follows since $\vct{w}_{\ell} = \pinv{(\matA_{\tau_\ell})} \z_{\ell}$. Therefore, we conclude that $\z_{\ell+1} = \vecb - \matA \x_{\ell+1}$ which completes the proof.
$\hfill\Box$

Combining this lemma with~\eqref{eq:zous2} yields the following result which shows convergence of the estimation to the least squares solution under the map $\mtx{A}$.
\begin{theorem}\label{thm:LSsimple}
Algorithm~\ref{alg:LSsimple} with input $\matA$, $\vecb$, $T\in\NN$, and $(\overline{p}, \overline{\alpha}, \overline{\beta})$ column paving $\mathcal{T}$, outputs an estimate vector $\x_{T}$ that satisfies
\[ \EE \norm{\matA (\xls - \x_{T})}^2 \leq \left(1 -\frac{\sigma_{\min}^2(\matA)}{ \overline{p}\overline{\beta}}\right)^{T} \norm{\br}^2.\]
\end{theorem}
\textit{Proof. }
 We observe that 
\begin{eqnarray*}
 \matA (\xls - \x_{(k)}) =   \br - \matA \x_{(k)} = \vecb - \matA \x_{(k)} - \bc  = \z_{(k)} - \bc
\end{eqnarray*}
where the first equality follows by $\br = \matA \pinv{\matA}\vecb = \matA\xls$, the second by orthogonality $\vecb = \bc + \br$ and the last equality from Lemma~\ref{claim:1}.  Combined with inequality~\eqref{eq:main} this yields the desired result.
$\hfill\Box$

When $\mtx{A}$ has full column rank, we may bound the estimation error $\|\xls - \x_T\|_2$ by $\frac{1}{\sigma_{\min}(\mtx{A})} \norm{\matA (\xls - \x_{T})}$ which combined with the fact that $\norm{\br}\leq \sigma_{\max}(\matA) \norm{\xls}$ implies the following corollary.

\begin{corollary}\label{cor:LSsimple}
Algorithm~\ref{alg:LSsimple} with full-rank $\matA$, $\vecb$, $T\in\NN$, and $(\overline{p}, \overline{\alpha}, \overline{\beta})$ column paving $\mathcal{T}$, outputs an estimate vector $\x_{T}$ that satisfies
\[ \EE \norm{ \xls - \x_{T}}^2 \leq \left(1 -\frac{\sigma_{\min}^2(\matA)}{ \overline{p}\overline{\beta}}\right)^{T} \kappa^2(\matA) \norm{\xls}^2.\]
\end{corollary}

\subsection{Implementation}\label{sec42}
The advantage to a single paving approach as in Algorithm~\ref{alg:LSsimple} is that one can utilize the
column-standardized version of $\matA$ while maintaining the same convergence to the (scaled version of the) least squares solution.  Utilizing Proposition~\ref{prop:intro-paving}, one is guaranteed a column-paving satisfying~\eqref{eq:paving}, so that paving parameters of Theorem~\ref{thm:LSsimple} and Corollary~\ref{cor:LSsimple} are bounded.  Since re-normalizing the columns of the matrix only re-scales the entries of $\xls$, Theorem~\ref{thm:LSsimple} and Corollary~\ref{cor:LSsimple} grant one access to the original least squares solution $\xls$.  To be precise, now let $\mtx{D}$ be the diagonal matrix whose entries correspond to the reciprocals of the column norms of $\matA$, so that $\overline{\matA} = \matA\mtx{D}$ has unit-norm columns.  Let ${\xlsn} = \overline{\matA}^\dagger\vecb$ denote the least squares solution of the re-normalized system, so that one has ${\xlsn} = \mtx{D}^{-1}\xls$.  Then
$$
\norm{\overline{\matA} ({\xlsn} - \x_{T})}^2 = \norm{\matA\mtx{D} (\mtx{D}^{-1}\xls - \x_{T})}^2 = \norm{\matA(\xls - \mtx{D}\x_{T})}^2.
$$
Thus applying Theorem~\ref{thm:LSsimple} for $\overline{\matA}$, and utilizing the fact that the range of $\matA$ is the same as that of $\overline{\matA}$, one has 
\[ \EE \norm{\matA(\xls - \mtx{D}\x_{T})}^2 \leq \left(1 -\frac{\sigma_{\min}^2(\overline{\matA})}{ \overline{p}\overline{\beta}}\right)^{T} \norm{\br}^2.\]
This then implies that when the matrix is full rank,
\begin{equation}\label{eq:mixed1}
 \EE \norm{ \xls - \mtx{D}\x_{T}}^2 \leq \left(1 -\frac{\sigma_{\min}^2(\overline{\matA})}{ \overline{p}\overline{\beta}}\right)^{T} \kappa^2(\matA) \norm{\xls}^2.
\end{equation}
Since $\overline{p}$ and $\overline{\beta}$ are the paving parameters of $\overline{A}$, one has by utilizing Proposition~\ref{prop:intro-paving} and substituting the bounds of~\eqref{eq:paving} into~\eqref{eq:mixed1} that
\begin{equation}\label{eq:mix}
 \EE \norm{ \xls - \mtx{D}\x_{T}}^2 \leq \left(1 -\frac{C'}{ \kappa^2(\overline{\matA})\log(1+d)}\right)^{T} \kappa^2(\matA)\norm{\xls}^2.
\end{equation}

Although this bound does depend on the conditioning of both $\matA$ and $\overline{\matA}$ (note the rate itself only depends on the conditioning of the latter), it is the first to guarantee linear convergence to the true least squares solution utilizing a paving which is guaranteed by Proposition~\ref{prop:intro-paving} while not placing any restrictions on the matrix $\matA$ itself (such as standardization).

\begin{remark}
In considering the improvements offered by both the REK method and the block Kaczmarz method, one may ask whether it is advantageous to run a traditional REK projection step as in~\eqref{eq:rek} along with a traditional block Kaczmarz update step as in~\eqref{eq:iterates}.  However, empirically we have observed that such a combination actually leads to a degradation in performance and requires far more epochs to converge than the algorithms discussed above.  We conjecture that it is important to run both the projection update and the Kaczmarz update ``at the same speed''; if the Kaczmarz update utilizes many rows at once, so should the projection update, and vice versa.
\end{remark}

\section{Summary of Approaches}\label{sec:sum}

Here we detail the various approaches proposed, and summarize the practical implementation in several frameworks. Since the block variants of the methods require the blocks be well conditioned, it is natural to rely on matrix pavings for the analysis.  Unfortunately, such pavings are only readily available when the matrix itself is properly normalized.  For that reason, we have proposed several practical alternatives for the important setting in which normalization is not feasible.  Both of our proposed methods still offer computational advantages over the standard approaches (see the next section).  We summarize these approaches here, which cover all possible settings.

\begin{description}
\item[Consistent systems:] When the system is consistent, one does not lose any convergence properties by re-scaling the system to be standardized, since the solution remains the same.  In this simple setting, one can utilize the standardized system (either standardizing a priori or on the fly), and benefit from the convergence guaranteed by Corollary~\ref{cor:LSrek}.
\item[Inconsistent standardized systems:] When the matrix is already standardized (as naturally occurs for example in Vandermonde matrices used in trigonometric approximation), as in the above case one can immediately utilize Corollary~\ref{cor:LSrek} to guarantee convergence to the least squares solution. 
\item[Inconsistent systems with bounded dynamic range:] If the system is not standardized, but the ratio of the largest to smallest row norm is bounded (as is the case in random matrices for example, which have tightly concentrated row norms), one can still utilize this result.  Indeed, by utilizing the same paving one would use if the matrix was actually standardized -- but not actually standardizing the system, Corollary~\ref{cor:LSrek} can be used to obtain the convergence rate given in \eqref{eq:dr}.  One sees that if the dynamic range $a_r$ is bounded (by say, a constant or even $\log n$), the method convergences to the original least squares solution with approximately the same convergence rate.
\item[Inconsistent systems with unbounded dynamic range:] If the matrix has a large variety of row norms, it may clearly be challenging to guarantee the desired row paving.  For that reason, it will be advantageous to use a column paving.  To that end, we propose Algorithm~\ref{alg:LSsimple} which is designed for systems that cannot be separated into well-conditioned row blocks.  As discussed in Section~\ref{sec42}, Theorem \ref{thm:LSsimple} can be used via the column-standardized paving to obtain the convergence rate given in \eqref{eq:mix}.  This bound guarantees linear convergence to the true least squares solution, even for matrices which are far from standardized.  The disadvantage of course is that the rate depends on the conditioning of the column-standardized version, which may be hard to explicitly bound (as is true for general large matrices anyway). 
\end{description}

In any of these cases, our results may be used to guarantee linear convergence in expectation to the true least squares solution of the system.  Moreover, because matrix blocks can be utilized, we often see a significant speedup in runtime due to practical considerations.  See the next section for examples of such behavior.

\section{Experimental Results}\label{sec:exper}

Here we present some experiments using simple examples to illustrate the benefits of block methods.  We do not claim optimized implementations of the method, and only run on small problem sizes; our purpose is only to demonstrate that even in these simple examples, the block method offers advantages to the standard method.  We refer the reader to~\cite{ZF12:REK,needell2013paved} for more empirical results for both REK and block methods.  

\begin{figure}[ht]
\begin{tabular}{cc}
\includegraphics[width=3in]{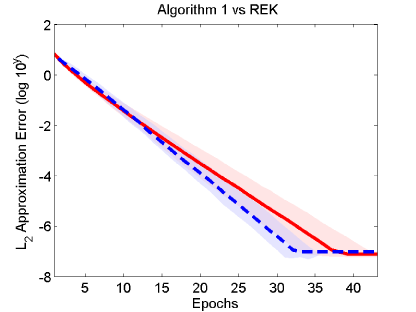} & \includegraphics[width=3in]{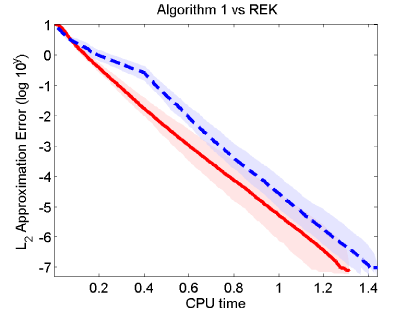} \\
\end{tabular}
\caption{$\ell_2$-norm error for REK (blue dashed) and Algorithm 1 (red) across epochs (left) and CPU time (right).  Matrix is $300\times 100$ Gaussian, system is consistent.\label{fig1}}
\end{figure}

\begin{figure}[ht]
\begin{tabular}{cc}
\includegraphics[width=3in]{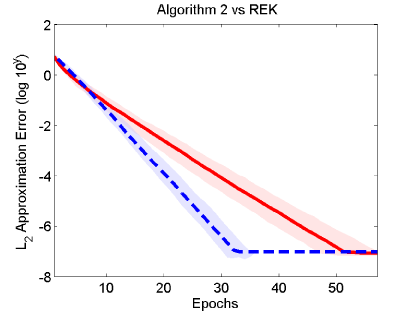} & \includegraphics[width=3in]{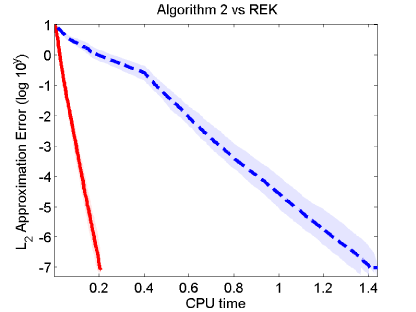} \\
\end{tabular}
\caption{$\ell_2$-norm error for REK (blue dashed) and Algorithm 2 (red) across epochs (left) and CPU time (right).  Matrix is $300\times 100$ Gaussian, system is consistent.\label{fig2}}
\end{figure}

In all experiments, one matrix is created and $40$ trials of each method are run.  In our first experiment, the matrix is a $300\times 100$ matrix with standard normal entries, whose rows are then normalized to each have norm one, yielding a condition number of $3.7$.  We use $30$ blocks, selected by a random partition.  The vector $\vct{x}$ is created to have independent standard normal entries, and the right hand side $\vct{b}$ is set to $\vct{Ax}$.  We track the $\ell_2$-error $\|\vct{x}_{LS}-\vct{x_k}\|_2$ across each epoch\footnote{We refer to an epoch as the number of iterations that is equivalent to one cycle through $n$ rows, even though rows and blocks are selected with replacement.  Thus for REK, an epoch is $n$ iterations, and for a block version with $b$ blocks, one epoch is $b$ iterations.  For an attempt at a fair comparison with Algorithm~\ref{alg:LSsimple} that only uses a column paving, we measure an ``epoch'' to be $n/b$ where $b$ is the number of blocks in the column paving.} as well as the CPU time (measured in Matlab using the cputime command).  In all experiments we considered a trial successful when the error reached $10^{-6}$.  The results for this case are presented in Figures~\ref{fig1} and \ref{fig2}.  In all figures, a heavy line represents median performance, and the shaded region spans the minimum to the maximum value across all trials.  As is demonstrated, even when the matrix does not have any natural block structure, the proposed algorithms outperform standard REK in terms of runtime.  

\begin{figure}[ht]
\begin{tabular}{cc}
\includegraphics[width=3in]{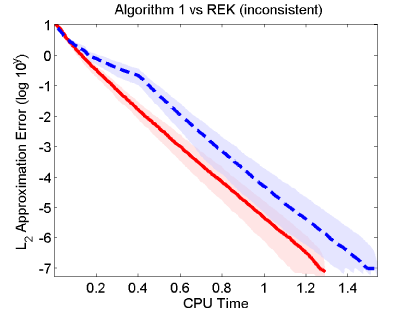} & \includegraphics[width=3in]{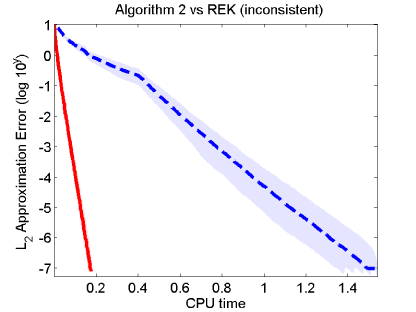} \\
\end{tabular}
\caption{Matrix is $300\times 100$ Gaussian, system is inconsistent. Left: $\ell_2$-norm error for REK (blue dashed) and Algorithm~\ref{alg:LSrek} (red) versus CPU time.  Right: $\ell_2$-norm error for REK (blue dashed) and Algorithm~\ref{alg:LSsimple} (red) versus CPU time.   \label{fig3}}
\end{figure}

\begin{figure}[ht]
\begin{tabular}{cc}
\includegraphics[width=3in]{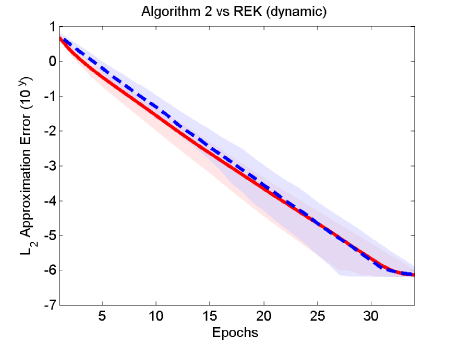} & \includegraphics[width=3in]{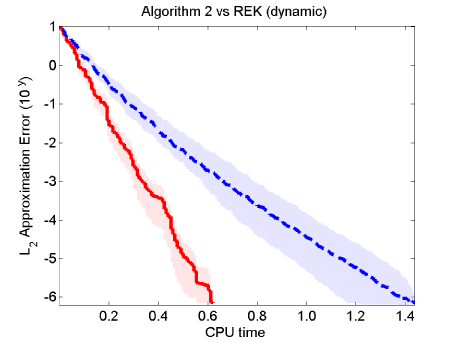} \\
\end{tabular}
\caption{Matrix is $300\times 100$ Gaussian with dynamic row norms, system is inconsistent. Plot shows $\ell_2$-norm error for REK (blue dashed) and Algorithm~\ref{alg:LSsimple} (red) versus epochs (left) and CPU time (right).   \label{figd}}
\end{figure}

Figure~\ref{fig3} shows similar plots, but in this case the system is no longer consistent.  For these experiments, we used the same type and size of the matrix $\mtx{A}$, but the right hand side vector $\vct{b}$ was generated as a Gaussian vector as well.  We created $\vct{b}$ so that the residual norm $\|\vct{b} - \vct{Ax_{LS}}\|_2 = 0.5$.  We then track the $\ell_2$-error between the iterate $\vct{x_k}$ and the least squares solution $\vct{x_{LS}}$ which we computed by $\mtx{A}^\dagger\vct{b}$.  We repeat the experiment also for a matrix whose dynamic range is not well bounded.  For that experiment, we generate a Gaussian matrix and then scale the row norms so that the $i$th row has norm equal to $i$.  The results for this case using Algorithm~\ref{alg:LSsimple} are show in Figure \ref{figd}.  The behavior in both cases, as predicted by our main results, is quite similar to the consistent case and thus breaks the convergence horizon of the standard Kaczmarz method.

\begin{figure}[ht]
\begin{tabular}{cc}
\includegraphics[width=3in]{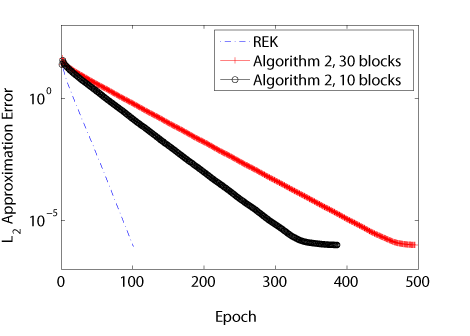} & \includegraphics[width=3in]{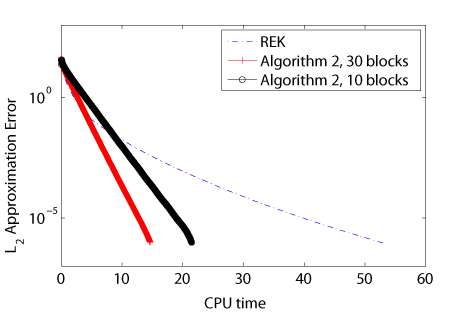} \\
\end{tabular}
\caption{System with $1200\times 400$ tomography matrix. Left: median $\ell_2$-norm error versus ``epoch'' (see footnote above).  Right: median $\ell_2$-norm error versus CPU time.   \label{fig4}}
\end{figure}

Lastly, we tested the methods on tomography problems, generated using the Matlab Regularization Toolbox by P.C. Hansen (\url{http://www.imm.dtu.dk/~pcha/Regutools/}) \cite{hansen2007regularization}.  In particular we present a 2D tomography problem $\mtx{A}\vct{x} = \vct{b}$ for an $n\times d$ matrix with $n=fN^2$ and $d=N^2$.  Here $\mtx{A}$ corresponds to the absorption along a random line through an $N\times N$ grid.  In our experiments we set $N=20$ and the oversampling factor $f=3$.  This yielded a matrix $\mtx{A}$ with condition number $\kappa(\mtx{A}) = 2.08$.  Since for this matrix it may be difficult to obtain a row paving, we instead use Algorithm~\ref{alg:LSsimple} and obtain a random column paving from the standardized version and use that for the matrix $\mtx{A}$.  The results for various choices of paving size are displayed in Figure~\ref{fig4}, which are in line with previous experiments.

\section{Related Work and Discussion}\label{sec:discuss}

The Kaczmarz method was first introduced in the 1937 work of Kaczmarz himself~\cite{K37:Angena}.  Since then, the method has been revitalized by researchers in computer tomography, under the name \textit{Algebraic Reconstruction Technique} (ART)~\cite{GBH70:Algebraic-Reconstruction,Byr08:Applied-Iterative,Nat01:Mathematics-Computerized,herman2009fundamentals}.  Deterministic convergence results for the method often depend on properties of the matrix that are difficult to compute or analyze~\cite{deutsch1985rate,deutsch1997rate,XZ02:Method-Alternating,galantai2005rate}.  Moreover, it has been well observed that random choice of row selection often speeds up the convergence~\cite{HS78:Angles-Null,HM93:Algebraic-Reconstruction,CFMSS92:New-Variants,Nat01:Mathematics-Computerized}. 

Recently, Strohmer and Vershynin~\cite{SV09:Randomized-Kaczmarz} derived the first provable convergence rate of the Kaczmarz method, showing that when each row is selected with probability proportional to its norm the method exhibits the expected linear convergence of~\eqref{eqn:RVrate}.  This work was extended to the inconsistent case in~\cite{Nee10:Randomized-Kaczmarz}, which shows linear convergence to within some fixed radius of the least squares solution.  The almost-sure guarantees were recently derived by Chen and Powell~\cite{CP12:Almost-Sure-Convergence}.  To break the convergence barrier, relaxation parameters can be introduced, so that each iterate is over or under projected onto each solution space.  Whitney and Meany prove that if the relaxation parameters tend to zero that the iterates converge to the least squares solution~\cite{whitney1967two}.  Further results using relaxation have also been obtained, see for example~\cite{censor1983strong,tanabe1971projection,hanke1990acceleration,ZF12:REK}.
An alternative to relaxation parameters was recently proposed by Zouzias and Freris~\cite{ZF12:REK} as the REK method described by~\eqref{eq:rek}.  Rather than alter the projection step, motivated by ideas of Popa~\cite{Pop98:Extensions-Block-Projections} they introduce a secondary step which aims to reduce the residual. 

The Kaczmarz method has been extended beyond linear systems as well.  For example, Leventhal and Lewis~\cite{LL10:Randomized-Methods} analyze the method for systems with polyhedral constraints and inequalities, which was also extended to the block case~\cite{briskman2014block}, and Richt{\'a}rik and Tak{\'a}{\v c}~\cite{RT11:Iteration-Complexity} build on these results for general optimization problems. 

Another important aspect of research in this area focuses on accelerating the convergence of the methods.  Geometric brute force methods can be used~\cite{EN11:Acceleration-Randomized}, additional row directions may be added~\cite{popa2012kaczmarz}, or instead one can select \textit{blocks} of rows rather than a single row in each iteration.  The block version of the Kaczmarz method is originally due to work of Elfving~\cite{Elf80:Block-Iterative-Methods} and Eggermont et al.~\cite{EHL81:Iterative-Algorithms}.  Its convergence rates were recently studied in~\cite{NW12:Two-Subspace-Projection} and analyzed via pavings by Needell and Tropp~\cite{needell2013paved}.  The block Kaczmarz method is of course a special instance in a broader class of block projection algorithms, see for example~\cite{XZ02:Method-Alternating} for a more general analysis and~\cite{Byr08:Applied-Iterative} for a presentation of other block variants.

To use block methods effectively, one needs to obtain a suitable partition of the rows (and/or columns).  Popa constructs such partitions by creating orthogonal blocks~\cite{Pop99:Block-Projections-Algorithms,Pop01:Fast-Kaczmarz-Kovarik,Pop04:Kaczmarz-Kovarik-Algorithm}, whereas Needell and Tropp promote the use of row pavings to construct the partition~\cite{needell2013paved}.  

Construction of pavings has been studied for quite some time now, and most early results rely on random selection.  The guarantee of lower and upper paving bounds has been derived by Bourgain and Tzafriri~\cite{BT87:Invertibility-Large} and Kashin and Tzafriri~\cite{KT94:Some-Remarks}, respectively.  Simultaneous guarantees were later derived by Bourgain and Tzafriri~\cite{BT91:Problem-Kadison-Singer} with suboptimal dependence on the matrix norm.  Recently, Spielman and Srivastava~\cite{SS12:Elementary-Proof} and Youssef~\cite{You12:Restricted-Invertibility} provided simple proofs of the results from~\cite{BT87:Invertibility-Large} and~\cite{KT94:Some-Remarks}, respectively.  Vershynin~\cite{Ver01:Johns-Decompositions} and Srivastava~\cite{Sri10:Spectral-Sparsification} extend the paving results to general matrices with arbitrary row norms; see also~\cite{You12:Restricted-Invertibility,You12:Note-Column}.  Proposition~\ref{prop:intro-paving} follows from the work of Vershynin~\cite{Ver06:Random-Sets} and Tropp~\cite{Tro09:Column-Subset}, and is attributed to the seminal work of Bourgain and Tzafriri~\cite{BT87:Invertibility-Large,BT91:Problem-Kadison-Singer}.  For particular classes of matrices, the paving can even be obtained from a random partition of the rows with high probability.  This is proved by Tropp~\cite{Tro08:Norms-Random} using ideas from~\cite{BT91:Problem-Kadison-Singer,Tro08:Random-Paving}, and is refined in~\cite{CD12:Invertibility-Submatrices}.

\section*{Acknowledgments}
D.N. is thankful to the Simons Foundation Collaboration grant and the Alfred P. Sloan Fellowship.  A.Z. has received funding from the European Research Council under the European Union's Seventh Framework Program (FP7/2007-2013) / ERC grant agreement $n^{o}$ 259569.  We would also like to thank Anna Ma for thoughtful discussions, and the reviewers for useful comments which significantly improved the manuscript.

\appendix\label{app}

\section{Proof of intermediate results}

Here we include the proof of Lemma~\ref{lem:elem}.  

\textit{Proof. }[Proof of Lemma~\ref{lem:elem}]

Assume the bounds~\eqref{eq:bnds} hold.  Applying the first bound in~\eqref{eq:bnds} recursively yields

\begin{align*}
\EE \|\vct{x}_{k^*} - \xls\|_2^2 &\leq \gamma^{k^*}\|\vct{x}_0 - \xls\|_2^2 + \sum_{j=0}^{k^*-1}\gamma^{k^*-1-j} r_j\\
&\leq \gamma^{k^*}\|\vct{x}_0 - \xls\|_2^2 + \sum_{j=0}^{\infty}\gamma^{j} B\\
&\leq \gamma^{k^*}\|\vct{x}_0 - \xls\|_2^2 + \frac{B}{1-\gamma},
\end{align*}

where the second inequality holds by the assumption that $r_k \leq \gamma^k B \leq B$, and the last by the properties of the geometric summation.  Similarly, observe that for any $k$ and $k^*$ we have

\begin{align*}
\EE \|\vct{x}_{k+k^*} - \xls\|_2^2 &\leq \gamma^{k}\EE\|\vct{x}_{k^*} - \xls\|_2^2 + \sum_{j=0}^{k-1}\gamma^{k-1-j} r_{j+k^*}\\
&\leq \gamma^{k}\EE\|\vct{x}_{k^*} - \xls\|_2^2 + \overline{\gamma}^{k^*}\sum_{j=0}^{\infty}\gamma^{j} B\\
&\leq \gamma^{k}\EE\|\vct{x}_{k^*} - \xls\|_2^2 + \overline{\gamma}^{k^*}\frac{B}{1-\gamma}.
\end{align*}

Now we choose $k$ and $k^*$ such that $T = k + k^*$ and $k=k^*$ if $T$ is even, or $k = k^*+1$ if $T$ is odd.  Combining the two inequalities above, we have

\begin{align*}
\EE \|\vct{x}_T - \xls\|_2^2 &= \EE \|\vct{x}_{k+k^*} - \xls\|_2^2\\
&\leq \gamma^{k}\EE\|\vct{x}_{k^*} - \xls\|_2^2 + \overline{\gamma}^{k^*}\frac{B}{1-\gamma}\\
&\leq \gamma^k\left( \gamma^{k^*}\|\vct{x}_0 - \xls\|_2^2 + \frac{B}{1-\gamma}\right) + \overline{\gamma}^{k^*}\frac{B}{1-\gamma}\\
&= \gamma^{k+k^*}\|\vct{x}_0 - \xls\|_2^2 + \left(\gamma^{k} + \overline{\gamma}^{k^*}\right)\frac{B}{1-\gamma}\\
&\leq \gamma^T\|\vct{x}_0 - \xls\|_2^2 + \left(\gamma^{\floor{T/2}} + \overline{\gamma}^{\floor{T/2}}\right)\frac{B}{1-\gamma}.
\end{align*}

This completes the proof.
$\hfill\Box$

%
%
\bibliographystyle{myalpha}
\bibliography{rk}

\newcommand{\etalchar}[1]{$^{#1}$}
\begin{thebibliography}{CFM{\etalchar{+}}92}

\bibitem[Bj{\"o}96]{Bjo96:Numerical-Methods}
{\AA}.~Bj{\"o}rck.
\newblock {\em Numerical Methods for Least Squares Problems}.
\newblock SIAM, Philadelphia, 1996.

\bibitem[BN]{briskman2014block}
J.~Briskman and D.~Needell.
\newblock Block kaczmarz method with inequalities.
\newblock {\em J. Math. Imaging Vis.}
\newblock To appear.

\bibitem[BT87]{BT87:Invertibility-Large}
J.~Bourgain and L.~Tzafriri.
\newblock Invertibility of ``large'' submatrices with applications to the
  geometry of {B}anach spaces and harmonic analysis.
\newblock {\em Israel J. Math.}, 57(2):137--224, 1987.

\bibitem[BT91]{BT91:Problem-Kadison-Singer}
J.~Bourgain and L.~Tzafriri.
\newblock On a problem of {K}adison and {S}inger.
\newblock {\em J. Reine Angew. Math.}, 420:1--43, 1991.

\bibitem[Byr08]{Byr08:Applied-Iterative}
C.~L. Byrne.
\newblock {\em Applied iterative methods}.
\newblock A K Peters Ltd., Wellesley, MA, 2008.

\bibitem[CD12]{CD12:Invertibility-Submatrices}
S.~Chr{\'e}tien and S.~Darses.
\newblock Invertibility of random submatrices via tail decoupling and a matrix
  {C}hernoff inequality.
\newblock {\em Statist. Probab. Lett.}, 82(7):1479--1487, 2012.

\bibitem[CEG83]{censor1983strong}
Y.~Censor, P.~P.~B. Eggermont, and D.~Gordon.
\newblock Strong underrelaxation in kaczmarz's method for inconsistent systems.
\newblock {\em Numer. Math.}, 41(1):83--92, 1983.

\bibitem[CFM{\etalchar{+}}92]{CFMSS92:New-Variants}
C.~Cenker, H.~G. Feichtinger, M.~Mayer, H.~Steier, and T.~Strohmer.
\newblock New variants of the {POCS} method using affine subspaces of finite
  codimension, with applications to irregular sampling.
\newblock In {\em Proc. SPIE: Visual Communications and Image Processing},
  pages 299--310, 1992.

\bibitem[CP12]{CP12:Almost-Sure-Convergence}
X.~Chen and A.~Powell.
\newblock Almost sure convergence of the {K}aczmarz algorithm with random
  measurements.
\newblock {\em J. Fourier Anal. Appl.}, pages 1--20, 2012.
\newblock 10.1007/s00041-012-9237-2.

\bibitem[Deu85]{deutsch1985rate}
F.~Deutsch.
\newblock Rate of convergence of the method of alternating projections.
\newblock {\em Parametric optimization and approximation}, 76:96--107, 1985.

\bibitem[DH97]{deutsch1997rate}
F.~Deutsch and H.~Hundal.
\newblock The rate of convergence for the method of alternating projections,
  ii.
\newblock {\em J. Math. Anal. Appl.}, 205(2):381--405, 1997.

\bibitem[EHL81]{EHL81:Iterative-Algorithms}
P.~P.~B. Eggermont, G.~T. Herman, and A.~Lent.
\newblock Iterative algorithms for large partitioned linear systems, with
  applications to image reconstruction.
\newblock {\em Linear Algebra Appl.}, 40:37--67, 1981.

\bibitem[Elf80]{Elf80:Block-Iterative-Methods}
T.~Elfving.
\newblock Block-iterative methods for consistent and inconsistent linear
  equations.
\newblock {\em Numer. Math.}, 35(1):1--12, 1980.

\bibitem[EN11]{EN11:Acceleration-Randomized}
Y.~C. Eldar and D.~Needell.
\newblock Acceleration of randomized {K}aczmarz method via the
  {J}ohnson-{L}indenstrauss lemma.
\newblock {\em Numer. Algorithms}, 58(2):163--177, 2011.

\bibitem[Gal05]{galantai2005rate}
A.~Gal{\'a}ntai.
\newblock On the rate of convergence of the alternating projection method in
  finite dimensional spaces.
\newblock {\em J. Math. Anal. Appl.}, 310(1):30--44, 2005.

\bibitem[GBH70]{GBH70:Algebraic-Reconstruction}
R.~Gordon, R.~Bender, and G.~T. Herman.
\newblock Algebraic reconstruction techniques ({ART}) for three-dimensional
  electron microscopy and {X}-ray photography.
\newblock {\em J. Theoret. Biol.}, 29:471--481, 1970.

\bibitem[Han07]{hansen2007regularization}
P.~C. Hansen.
\newblock Regularization tools version 4.0 for matlab 7.3.
\newblock {\em Numer. Algorithms}, 46(2):189--194, 2007.

\bibitem[Her09]{herman2009fundamentals}
G.~T. Herman.
\newblock {\em Fundamentals of computerized tomography: image reconstruction
  from projections}.
\newblock Springer, 2009.

\bibitem[HM93]{HM93:Algebraic-Reconstruction}
G.~Herman and L.~Meyer.
\newblock Algebraic reconstruction techniques can be made computationally
  efficient.
\newblock {\em IEEE Trans. Medical Imaging}, 12(3):600--609, 1993.

\bibitem[HN90]{hanke1990acceleration}
M.~Hanke and W.~Niethammer.
\newblock On the acceleration of kaczmarz's method for inconsistent linear
  systems.
\newblock {\em Linear Algebra Appl.}, 130:83--98, 1990.

\bibitem[HS78]{HS78:Angles-Null}
C.~Hamaker and D.~C. Solmon.
\newblock The angles between the null spaces of {X}-rays.
\newblock {\em J. Math. Anal. Appl.}, 62(1):1--23, 1978.

\bibitem[Kac37]{K37:Angena}
S.~Kaczmarz.
\newblock Angen\"aherte aufl\"osung von systemen linearer gleichungen.
\newblock {\em Bull. Internat. Acad. Polon.Sci. Lettres A}, pages 335--357,
  1937.

\bibitem[KT94]{KT94:Some-Remarks}
B.~Kashin and L.~Tzafriri.
\newblock Some remarks on coordinate restriction of operators to coordinate
  subspaces.
\newblock Insitute of Mathematics Preprint~12, Hebrew University, Jerusalem,
  1993--1994.

\bibitem[LL10]{LL10:Randomized-Methods}
D.~Leventhal and A.~S. Lewis.
\newblock Randomized methods for linear constraints: convergence rates and
  conditioning.
\newblock {\em Math. Oper. Res.}, 35(3):641--654, 2010.

\bibitem[LWS14]{liu2014asynchronous}
J.~Liu, S.~J. Wright, and S.~Sridhar.
\newblock An asynchronous parallel randomized kaczmarz algorithm.
\newblock {\em arXiv preprint arXiv:1401.4780}, 2014.

\bibitem[Nao11]{Nao11:Sparse-Quadratic}
A.~Naor.
\newblock Sparse quadratic forms and their geometric applications.
\newblock Technical Report No.~1033, S{\'e}minaire Bourbaki, Jan. 2011.

\bibitem[Nat01]{Nat01:Mathematics-Computerized}
F.~Natterer.
\newblock {\em The mathematics of computerized tomography}, volume~32 of {\em
  Classics in Applied Mathematics}.
\newblock Society for Industrial and Applied Mathematics (SIAM), Philadelphia,
  PA, 2001.
\newblock Reprint of the 1986 original.

\bibitem[Nee10]{Nee10:Randomized-Kaczmarz}
D.~Needell.
\newblock Randomized {K}aczmarz solver for noisy linear systems.
\newblock {\em BIT}, 50(2):395--403, 2010.

\bibitem[NT13]{needell2013paved}
D.~Needell and J.~A. Tropp.
\newblock Paved with good intentions: Analysis of a randomized block kaczmarz
  method.
\newblock {\em Linear Algebra Appl.}, pages 199--221, 2013.

\bibitem[NW13]{NW12:Two-Subspace-Projection}
D.~Needell and R.~Ward.
\newblock Two-subspace projection method for coherent overdetermined linear
  systems.
\newblock {\em J. Fourier Anal. Appl.}, 19(2):256--269, 2013.

\bibitem[Pop98]{Pop98:Extensions-Block-Projections}
C.~Popa.
\newblock Extensions of block-projections methods with relaxation parameters to
  inconsistent and rank-deficient least-squares problems.
\newblock {\em BIT}, 38(1):151--176, 1998.

\bibitem[Pop99]{Pop99:Block-Projections-Algorithms}
C.~Popa.
\newblock Block-projections algorithms with blocks containing mutually
  orthogonal rows and columns.
\newblock {\em BIT}, 39(2):323--338, 1999.

\bibitem[Pop01]{Pop01:Fast-Kaczmarz-Kovarik}
C.~Popa.
\newblock A fast {K}aczmarz-{K}ovarik algorithm for consistent least-squares
  problems.
\newblock {\em Korean J. Comput. Appl. Math.}, 8(1):9--26, 2001.

\bibitem[Pop04]{Pop04:Kaczmarz-Kovarik-Algorithm}
C.~Popa.
\newblock A {K}aczmarz-{K}ovarik algorithm for symmetric ill-conditioned
  matrices.
\newblock {\em An. \c Stiin\c t. Univ. Ovidius Constan\c ta Ser. Mat.},
  12(2):135--146, 2004.

\bibitem[PPKR12]{popa2012kaczmarz}
C.~Popa, T.~Preclik, H.~K{\"o}stler, and U.~R{\"u}de.
\newblock On {K}aczmarz's projection iteration as a direct solver for linear
  least squares problems.
\newblock {\em Linear Algebra Appl.}, 436(2):389--404, 2012.

\bibitem[RT11]{RT11:Iteration-Complexity}
P.~Richt{\'a}rik and M.~Tak{\'a}{\v c}.
\newblock Iteration complexity of randomized block-coordinate descent methods
  for minimizing a composite function.
\newblock Available at \url{arXiv:1107.2848}, Apr. 2011.

\bibitem[Sri10]{Sri10:Spectral-Sparsification}
N.~Srivastava.
\newblock {\em Spectral sparsification and restricted invertibility}.
\newblock Phd dissertation, Yale University, New Haven, CT, 2010.

\bibitem[SS87]{SS87:Applications}
K.~M. Sezan and H.~Stark.
\newblock Applications of convex projection theory to image recovery in
  tomography and related areas.
\newblock In H.~Stark, editor, {\em Image Recovery: {T}heory and application},
  pages 415–--462. Acad. Press, 1987.

\bibitem[SS12]{SS12:Elementary-Proof}
D.~A. Spielman and N.~Srivastava.
\newblock An elementary proof of the restricted invertibility theorem.
\newblock {\em Israel J. Math.}, 190:83--91, 2012.

\bibitem[SV06]{SV06:Arandom}
T.~Strohmer and R.~Vershynin.
\newblock A randomized solver for linear systems with exponential convergence.
\newblock In {\em RANDOM 2006 (10th International Workshop on Randomization and
  Computation)}, number 4110 in Lecture Notes in Computer Science, pages
  499--507. Springer, 2006.

\bibitem[SV09]{SV09:Randomized-Kaczmarz}
T.~Strohmer and R.~Vershynin.
\newblock A randomized {K}aczmarz algorithm with exponential convergence.
\newblock {\em J. Fourier Anal. Appl.}, 15(2):262--278, 2009.

\bibitem[Tan71]{tanabe1971projection}
K.~Tanabe.
\newblock Projection method for solving a singular system of linear equations
  and its applications.
\newblock {\em Numer. Math.}, 17(3):203--214, 1971.

\bibitem[Tro08a]{Tro08:Norms-Random}
J.~A. Tropp.
\newblock Norms of random submatrices and sparse approximation.
\newblock {\em C. R. Math. Acad. Sci. Paris}, 346(23-24):1271--1274, 2008.

\bibitem[Tro08b]{Tro08:Random-Paving}
J.~A. Tropp.
\newblock The random paving property for uniformly bounded matrices.
\newblock {\em Studia Math.}, 185(1):67--82, 2008.

\bibitem[Tro09a]{Tro09:Column-Subset}
J.~A. Tropp.
\newblock Column subset selection, matrix factorization, and eigenvalue
  optimization.
\newblock In {\em Proceedings of the {T}wentieth {A}nnual {ACM}-{SIAM}
  {S}ymposium on {D}iscrete {A}lgorithms}, pages 978--986, Philadelphia, PA,
  2009. SIAM.

\bibitem[Tro09b]{tropp2009column}
J.~Tropp.
\newblock Column subset selection, matrix factorization, and eigenvalue
  optimization.
\newblock In {\em Proceedings of the twentieth Annual ACM-SIAM Symposium on
  Discrete Algorithms}, pages 978--986. Society for Industrial and Applied
  Mathematics, 2009.

\bibitem[Ver01]{Ver01:Johns-Decompositions}
R.~Vershynin.
\newblock John's decompositions: selecting a large part.
\newblock {\em Israel J. Math.}, 122:253--277, 2001.

\bibitem[Ver06]{Ver06:Random-Sets}
R.~Vershynin.
\newblock Random sets of isomorphism of linear operators on {H}ilbert space.
\newblock In {\em High dimensional probability}, volume~51 of {\em IMS Lecture
  Notes Monogr. Ser.}, pages 148--154. Inst. Math. Statist., Beachwood, OH,
  2006.

\bibitem[WM67]{whitney1967two}
T.~M. Whitney and R.~K. Meany.
\newblock Two algorithms related to the method of steepest descent.
\newblock {\em SIAM J. Numer. Anal.}, 4(1):109--118, 1967.

\bibitem[XZ02]{XZ02:Method-Alternating}
J.~Xu and L.~Zikatanov.
\newblock The method of alternating projections and the method of subspace
  corrections in {H}ilbert space.
\newblock {\em J. Amer. Math. Soc.}, 15(3):573--597, 2002.

\bibitem[You12a]{You12:Note-Column}
P.~Youssef.
\newblock A note on column subset selection.
\newblock Available at \url{arXiv:1212.0976}, Dec. 2012.

\bibitem[You12b]{You12:Restricted-Invertibility}
P.~Youssef.
\newblock Restricted invertibility and the {B}anach--{M}azur distance to the
  cube.
\newblock Available at \url{arXiv:1206.0654}, June 2012.

\bibitem[ZF12]{ZF12:REK}
A.~Zouzias and N.~M. Freris.
\newblock Randomized extended {K}aczmarz for solving least-squares.
\newblock {\em SIAM J. Matrix Anal. A.}, 34(2):773--793, 2012.

\end{thebibliography}
%
%

\end{document}